\documentclass[11pt]{article}
\usepackage{geometry} 
\usepackage{graphicx}
\usepackage{amssymb}
\usepackage{epstopdf}
\usepackage{amsmath}
\usepackage{amsthm}
\usepackage{fullpage}

\title{On some classical constructions \\extended to hyperbolic geometry
}
\author{Arseniy V. Akopyan\thanks{Institute for Information Transmission Problems, Moscow, Russia. \textit{e-mail:}\texttt{akopjan@gmail.com}}}
\date{} 

\DeclareMathOperator{\angb}{{}_\angle}
\DeclareMathOperator{\san}{\sphericalangle}
\DeclareMathOperator{\trib}{{}_\triangle}

\newcounter{ris}
\renewcommand{\r}{\refstepcounter{ris}%
 Fig. \arabic{ris}}

\begin{document}
	\ifpdf
	\DeclareGraphicsExtensions{.pdf, .jpg, .tif, .mps}
	\else
	\DeclareGraphicsExtensions{.eps, .jpg, .mps}
	\fi	
	
\maketitle
\theoremstyle{theorem}
\newtheorem{theorem}{Theorem}
\newtheorem*{conjecture*}{Conjecture}
\newtheorem{lemma}[theorem]{Lemma}

\begin{abstract}
We consider some constructions in hyperbolic geometry that are analogous to classical constructions in Euclidean geometry. We show that both Monge's theorem and the theorem on the concurrence of the common chords of three circles also hold in absolute geometry. We proffer analogues of the Euler line and the Euler nine-point circle and also extend Feuerbach's famous theorem.
\end{abstract}

\section{Introduction}
In the past decade there has been a growing interest in so-called classical geometry. Many new and interesting theorems have been proved, and numerous books and articles have been written. However, there have been almost no investigations of even the most rudimentary (by today's standards) constructions in hyperbolic and spherical geometry.

In this article we formulate some non-Euclidean analogues of classical theorems of Euclidean geometry. The principal result is an analogue of Feuerbach's theorem regarding the tangency of the inscribed circle and the Euler nine-point circle (recall that this circle passes through the midpoints of the sides, the feet of the altitudes, and also the midpoints of the segments joining the vertices to the orthocenter).

In the next section, we formulate an analogue of this theorem. 
In section~\ref{sct:(pseudovisoti)} we explain the meaning of this analogue. 
In section~\ref{sct:(parallelnost)} we prove the existence of the Euler circle. 
Next, in section~\ref{sct:(stepentochki)}, we construct analogues of the concepts of the power of a point, inversion and dilation. 
We use these to prove the existence of the Euler line in section~\ref{sct:(eulerline)}. 
In section~\ref{sct:(feuerbach theorem)} we prove Feuerbach's theorem. 
In section~\ref{sct:(monge theorem)} we prove Monge's theorem and use it to prove an interesting property of the Feuerbach point.

Please note that we shall assume below the facts that any two lines intersect and that three points determine a circle. In hyperbolic geometry this is not generally true. Straight lines may diverge, and circles may degenerate into horocycles or equidistant curves (which may be considered circles with centers on the absolute or beyond respectively). Consideration of all possible cases would not only complicate the proof, but would contain no fundamentally new ideas. To complete our arguments, we could always say that other cases follow from a theorem by analytic continuation, since the cases considered by us are sufficiently general (they include an interior point in the configuration space).

Nevertheless, in the course of our argument we shall try to avoid major errors and show that the statements can be demonstrated without resorting to more powerful tools.

\section{Formulation}
\label{sct:(formulirovka)}

In 1822, Karl Wilhelm Feuerbach proved that the Euler nine-point circle was tangent to the inscribed and the three escribed circles of a triangle. This theorem is now known as Feuerbach's theorem. It has several interesting and varied proofs. For example, it follows directly from Casey's theorem, a generalization of Ptolemy's theorem. Another proof can be obtained through an inversion centered at the midpoint of a side of the triangle.

Vladimir Protasov in \cite{protasov1999thefeuerbach} provided a proof of Feuerbach's theorem that is not trivial, but elementary and very informative. He showed that it is a special case of the so-called \textit{segment theorem}.

There are many interesting generalizations of Feuerbach's theorem. See, e.g. \cite{emelyanov2002family_en} or \cite{markelov1999parabola_en}.

In 2008, P. V. Bibikov formulated a hypothesis that is a generalization of the theorem for hyperbolic geometry.

\begin{conjecture*}[Pavel Bibikov]
Let $M_a$, $M_b$ and $M_c$ be points on the sides $a$, $b$ and $c$ respectively of triangle $ABC$ such that the cevians $AM_a$, $BM_b$ and $CM_c$ bisect the area of the triangle. Then the circle passing through $M_a$, $M_b$ and $M_c$ is tangent to the inscribed and three escribed circles of the triangle.

\end{conjecture*}

We call this circle the {\em Euler circle} of the triangle. We henceforth refer to the area-bisecting cevians $AM_a$, $BM_b$ and $CM_c$ as the {\em bisectors} of the triangle.

We shall prove this conjecture and in addition show the collinearity of (1) the circumcenter of a triangle, (2) the intersection of its bisectors, (3) the center of the Euler circle and (4) the intersection of the cevians through the remaining intersections of the Euler circle with the sides of the triangle. The line containing these four points, in the opinion of the author, deserves to be called the {\em Euler line} of the triangle.

\section{Pseudomedians and pseudoaltitudes}
\label{sct:(pseudovisoti)}

Before we deal with Feuerbach's theorem, we consider the six points on the sides of a triangle where it is intersected by the Euler circle. Three of these intersections are on the bisectors, which divide the area of the triangle in half. In Euclidean geometry these lines are also the medians, so we could refer to the bisectors as {\em pseudomedians}. In hyperbolic geometry, there are identical angle sums for the two triangles into which each bisector divides the triangle. Accordingly it is natural to formulate some definition of the altitudes of a Euclidean triangle in terms of these angles and generalize the definition for hyperbolic geometry.

Let $H_a$ be the foot of the altitude from vertex $A$ to side $BC$. Then, in the Euclidean case, we have the following:
\begin{equation*}
\begin{split}
\angb AH_aB - (\angb H_aAB+\angb ABH_a)=
\angb AH_aC - (\angb H_aAC+\angb ACH_a) \\(= {1\over2}(\pi - (\angb A+\angb B +\angb C))=0)\text{.}
\end{split}
\end{equation*}

Call cevians with a similar property in hyperbolic geometry {\em pseudoaltitudes}. 
Let us verify the existence of pseudoaltitudes.

Denote by $\san XYZ$ the value of $\angb XYZ-\angb ZXY-\angb YZX$, where the angles are signed; \mbox{$\angb XYZ>0$} if and only if the rotation around $Y$ along the shortest arc from $X$ to $Z$ is counterclockwise.

Consider triangle $ABC$, in which the order of the vertices $A$, $B$, $C$ is clockwise. Using the above rules for signs, it is possible to verify that when moving point $X$ on line $BC$ in the direction from $B$ to $C$, the value of $\san BXA$ is continuously decreasing, while the value of $\san AXC$ is continuously increasing; $\san BXA-\san AXC>0$ when $X$ is far enough out on the ray $CB$, and \mbox{$\san BXA-\san AXC<0$} when $X$ is far enough out on the ray $BC$. Therefore there must be a point~$H_a$ for which $\san BXA=\san AXC$, and this point is the foot of the pseudoaltitude from vertex~$A$.

Let $H_b$ and $H_c$ be the feet of the pseudoaltitudes on lines $CA$ and $AB$ respectively. This requires that $\san CH_bB=\san BH_bA$ and $\san AH_cC=\san CH_cB$. In the next section, we show that the Euler circle passes through these three points.

Note that
\begin{equation*}
 \san BH_aA+\san AH_aC= \san
CH_bB+\san BH_bA =\san AH_cC+\san
CH_cB=S_{\trib ABC},
\end{equation*}
from which it follows that
\begin{equation*}
 \san BH_aA=\san AH_aC= \san
CH_bB=\san BH_bA= \san AH_cC=\san
CH_cB={1\over2}S_{\trib ABC}.
\end{equation*}

\begin{center}
 \includegraphics{figfeuerbach-1}

 \r
\end{center}

It is easy to prove in absolute geometry that quadrilaterals $ABH_aH_b$, $ACH_aH_c$ and $BCH_bH_c$ are cyclic. (See F. V. Petrov \cite{petrov2009cyclic} for a fuller discussion of cyclic quadrilaterals.) To do this we use an analogue of the theorem of the inscribed angle. This theorem is in many books (see, for example, \cite{scopets1962exercises_en} or \cite{prasolov2001geometry} or the aforementioned article \cite{petrov2009cyclic}), but we nevertheless recall its proof.

\begin{lemma}
	For any point $X$ on arc $AB$ of a fixed circle, the value of $\san AXB$ is constant.
\end{lemma}
\begin{proof}
	Let $O$ be the center of a circle and $X$ be an arbitrary point on the arc $AB$. Obviously any line from $O$ to $AB$ is perpendicular to $AB$, and in addition, we have
	\begin{equation*}
	\angb OAB = {1\over2}(\angb AXB -\angb BAX-\angb
	XBA)={1\over2}\san AXB\text{.}
	\end{equation*}
	So, the value of $\san AXB$ does not depend on the choice of $X$ and is always equal to $2\angb OAB$.
\end{proof}

\begin{center}
 \includegraphics{figfeuerbach-2}

 \r
\end{center}

We can also show the converse, that if $\san AXB=\san AYB$, then the four points $A$, $B$, $X$ and $Y$ are concyclic in the Poincar\'e model, i.e., their projections lie on a straight line, a circle, a horocycle, or an equidistant curve.

What will benefit us is a consequence of the inscribed angle theorem.

\begin{lemma}
	There exists a unique triple of points $X_a$, $X_b$ and $X_c$ on the sides of triangle $ABC$ such that the quadrilaterals $ABX_aX_b$, $ACX_cX_a$ and $BCX_cX_b$ are cyclic.
\label{lem:(edinstvennost)}
\end{lemma}
\begin{proof}
	Indeed, from the inscribed angle theorem we obtain
\begin{equation*}
	\begin{split}
	 \san BX_aA=\san BX_bA= S_{\trib ABC} -
	\san CX_bB= S_{\trib ABC} - \san CX_cB=\\
	=\san AX_cC=\san AX_aC=S_{\trib ABC} -
	\san BX_aA.
	\end{split}
\end{equation*}

	That is, $\san BX_aA={1\over2}S_{\trib ABC}=\san BH_aA$. This means that points $X_a$ and $H_a$ coincide. A~similar argument applies to the other two points.
\end{proof}

\section{Parallelism and antiparallelism}
\label{sct:(parallelnost)}

Recall that, with respect to a Euclidean angle $XOY$, two lines can be antiparallel. 
Let the line $OX$ include points $A$ and $B$, the line $OY$ include the points $C$ and $D$, and the quadrilateral $ABCD$ be cyclic.
In such a case we say that the segments $AD$ and $BC$ are antiparallel.
In Euclidean geometry, a quadrilateral is cyclic if and only if the sum of its opposite angles is equal to $\pi$. Thus, it is easy to see that if two segments $A'D'$ and $B'C'$ are parallel to two antiparallel segments $AD$ and $BC$ respectively, then they will also be antiparallel to each other.

\begin{center}
 \includegraphics{figfeuerbach-3}

 \r
\end{center}

We carry this concept over to the hyperbolic plane. Just as in the Euclidean case, let the segments $AD$ and $CD$ be antiparallel and the quadrilateral $ABCD$ be cyclic. Analogous to the Euclidean case, we can use the criterion of concyclicity and find that the segments are antiparallel if and only if $\angb DAO-\angb ODA=\angb OCB-\angb CBO$, which is equivalent to $\san DAO=\san OCB$.

We can now introduce the concept of parallel segments relative to an angle. The segments $BC$ and $AD$ are parallel if
\begin{equation*}
	\angb DAO-\angb ODA=\angb CBO-\angb OCB,
\end{equation*}
which is equivalent to
\[
	\san DAO=-\san OCB.
\]
It is clear that any two segments that are antiparallel to some segment must be parallel to each other.

Note that in the article by F. V. Petrov \cite{petrov2009cyclic}, a quadrilateral $ABCD$ with parallel segments $AD$ and $BC$ is called a trapezoid.  The following theorem about it has been proven.
\begin{lemma}[{\cite{petrov2009cyclic}, Lemma 3}]
	\label{lem:trapezia}
	The area of triangles $ABC$ and $ABD$ are equal if and only if the following condition holds for the sum of the angles of quadrilateral $ABCD$:
\begin{equation*}
	\angb A+\angb D=\angb B+\angb C.
\end{equation*}
\end{lemma}

From this lemma we see that when $M_a$ and $M_b$ are the feet of the bisectors, then the segments $M_aM_b$ and $AB$ are parallel with respect to angle $ACB$. Since $H_aH_b$ is antiparallel to $AB$ relative to the same angle, then $H_aH_b$ is antiparallel to $M_aM_b$.

Now it is easy to prove that $H_a$, $H_b$ and $H_c$ are also points of intersection of the circumcircle of $M_aM_bM_c$ with the sides of the triangle. Indeed, let $X_a$, $X_b$ and $X_c$ be the second points of intersection with the sides. Then $X_aX_b$ is antiparallel to $M_aM_b$, which means that it is also antiparallel to $AB$, i.e., the quadrilateral $X_aX_bAB$ is cyclic.   Similarly quadrilaterals $ACX_aX_b$ and $BCX_bX_c$ are cyclic. But by Lemma~\ref{lem:(edinstvennost)}, such a triple of points is unique, and must therefore be the same as the triple $H_a$, $H_b$ and $H_c$.

In the Poincar\'e disk there is natural connection between parallelism and antiparallelism. 
The latter, as we know, corresponds to the case of an inscribed quadrilateral.
A proof of the following theorem can be found for example in \cite{shvartsman2007comment_en}.

\begin{theorem}[analogue in hyperbolic geometry of Lexell's theorem]
	\label{thm:parallel}
Let points $A^*$ and $B^*$ be the inverses of points $A$ and $B$ with respect to the absolute circle of the Poincar\'e disk model. Let $\omega$ be any Euclidean circle passing through $A^*$ and $B^*$. Then for any point $X$ on $\omega$, the area of triangle $XAB$ (in the hyperbolic sense) is constant, i.e., will not depend on the choice of $X$.
\end{theorem}
It follows from this theorem and Lemma~\ref{lem:trapezia} that when segments $AB$ and $DC$ are parallel, then points $A^*$, $B^*$, $C$ and $D$ are concyclic in the Poincar\'e disk.

Note also that Lemma~\ref{lem:trapezia} can be deduced from Theorem~\ref{thm:parallel}.

\section{The power of a point, homothety and inversion}
\label{sct:(stepentochki)}

Let the {\em pseudolength} $d_E(A,B)$ of a hyperbolic line segment $AB$ be defined as follows. After a suitable motion that places $A$ at the center of a Poincar\'e disk of radius 1, the pseudolength $d_E(A,B)$ is defined as the Euclidean length of segment $AB$. 

It is easy to verify that the pseudolength does not depend on the motion that takes $A$ to the center of the Poincar\'e disk. Also, the pseudolength is symmetric: $d_E(A,B)=d_E(B,A)$.

With the help of pseudolengths, we can define the concepts of the power of a point with respect to a circle, homothety and inversion. Thus, the power of a point $A$ with respect to a circle $\omega$ is defined as the product of the pseudolengths of the segments $AB$ and $AC$, where $BC$ is a chord of the circle $\omega$ that passes through $A$. That this is well defined (i.e., the value of $d_E(A,B)\cdot d_E(A,C)$ is independent of the choice of the chord through $A$) becomes apparent after moving point $A$ to the center of the Poincar\'e disk, in which case the theorem is equivalent to a Euclidean theorem.

We can now define the radical axis of two circles as the set of points for which the powers of each point with respect to the two circles are the same.

\begin{lemma}
	The radical axis is a straight line.
	\label{lem:no8}
\end{lemma}
This is evident for two intersecting circles, as they have a common chord. For circles that do not intersect, we can consider them to be the intersection of a plane with two intersecting spheres in space. The radical axis of the circles is the intersection with this plane of the plane of the circle of intersection of the two spheres. For more on this construction, see \cite{akopyan2007geometry}.

We note that it is now possible to transfer to hyperbolic geometry any theorems of Euclidean geometry that depend upon the radical axis. For example, we can prove Brianchon's theorem for a circle in hyperbolic geometry. (See the proof of Brianchon's theorem employing the radical axis in \cite{akopyan2007geometry} or \cite{prasolov2006problems_en}.)

From Lemma~\ref{lem:no8} it follows immediately that the pseudoaltitudes are coincident at the radical center of circles $ABH_aH_b$, $ACH_aH_c$ and $BCH_bH_c$.

Similarly, the bisectors are coincident at the radical center of circles $A^*B^*M_aM_b$, $A^*C^*M_aM_c$ and $B^*C^*M_bM_c$ (because all of the lines through $A^*$ also pass through $A$, so that the line $A^*M_a$ is nothing but the bisector $AM_a$).

We now define a {\em homothety with center $A$ and ratio $k$} as the transformation taking point $P$ to point $P'$ lying on the ray $AP$ (in the case of positive $k$) such that $d_E(A,P')=k\cdot d_E(A,P)$.

Similarly we define an {\em inversion with center at $A$ and radius $r$} as the transformation that takes point $P$ to the point $P'$ lying on the ray $AP$ such that $d_E(A,P')\cdot d_E(A,P)=r^2$.

It is easy to see that these transformations take a circle in the Poincar\'e disk into a circle and a straight line through the center into a straight line.

\section{The Euler line}
\label{sct:(eulerline)}

We denote by $O$ the circumcenter of triangle $ABC$, by $E$ the center of its Euler circle, and by $M$ and $H$ the respective points of intersection of the bisectors and the pseudoaltitudes. We shall show that these points lie on a straight line, which we shall call the Euler line.

Move point $M$ to the center of the Poincar\'e disk. Then the criterion of cyclicity for hyperbolic and Euclidean quadrilaterals is that (Euclidean) straight lines $AB$ and $M_aM_b$ are parallel in the Euclidean sense. Consequently, we have the following equation:
\begin{equation*}
	\frac{d_E(M,M_a)}{d_E(M,A)}=\frac{d_E(M,M_b)}{d_E(M,B)}=\frac{d_E(M,M_c)}{d_E(M,C)}.
\end{equation*}
That is, the circumcircle of $ABC$ is transformed into a circle by the hyperbolic Euler homothety with center $M$ and ratio $-d_E(M,M_a)/d_E(M,A)$. Of course, under this map the center of the circumscribed circle is not mapped to the point $E$ (i.e., the center of the image circle), but nevertheless it must lie on $ME$. In addition, it must lie on the line $MO$. We see that these two lines coincide, so that the points $O$, $E$ and $M$ are collinear.

We employ a similar argument for the point $H$. Since $H$ is the radical center of circles $ABH_aH_b$, $ACH_aH_c$ and $BCH_bH_c$, it follows that
\begin{equation*}
	d_E(H,H_a)\cdot d_E(H,A)=d_E(H,H_b)\cdot d_E(H,B)=d_E(H,H_c)\cdot d_E(H,C).
\end{equation*}
That is, the Euler circle is the image of the circumcircle of triangle $ABC$ after rotations around the point $H$ through an angle of $\pi$ and an inversion with center $H$ and radius squared of $d_E(H,H_a)\cdot d_E(H,A)$. This implies that $H$ also lies on the line $OE$.

This argument is similar to the proof of the existence of the Euler line in Euclidean geometry. Therefore, in the opinion of the author, this line can be considered a direct analogue in hyperbolic geometry of the Euler line.

\section{Feuerbach's theorem}
\label{sct:(feuerbach theorem)}

We now state the segment theorem (see \cite{protasov1999thefeuerbach}), which will be the principal tool in proving Feuerbach's theorem.
\begin{theorem}[Segment theorem, Protasov \cite{protasov1999thefeuerbach}]
Given two lines $a$ and $b$, each tangent to a given circle $\omega$, let a third line $c$ be tangent to $\omega$ and intersect $a$ at $A$ and $b$ at $B$. Let $\omega_c$ be a circle through $A$ and $B$ such that the central angle of arc $AB$ is $\phi$. There are two circles that are each tangent to $\omega_c$ as well as $a$ and $b$ and remain fixed even as $c$ is varied.
\end{theorem}

Consider Feuerbach's theorem in the Poincar\'e disk model from a Euclidean point of view. Ignore the feet of the bisectors, which we do not need, and consider only the pseudoaltitudes $AA'$, $BB'$ and $CC'$ of triangle $ABC$. We will see the following theorem.

\begin{theorem}
	We are given two triangles $ABC$ and $A'B'C'$, such that quadrilaterals $ABB'A'$, $ACC'A'$ and $BCC'B'$ are cyclic. Then there are four circles that are tangent to the circumcircles of triangles $ABC'$, $AB'C$, $A'BC$ and $A'B'C'$.
\end{theorem}

\begin{center}
 \includegraphics{figfeuerbach-4}

 \r
\end{center}

\begin{proof}
	We make an inversion with center at $A$. Then the circles $ABC'$ and $AB'C$ are transformed into lines, and the rest of the circles remain circles. Thus, two of the four circles are now lines.

	So suppose there is the following construction. Points $B$, $B'$, $C'$ and $C$ lie on one circle. Lines $BB'$ and $CC'$ intersect at $A'$. We need to prove that there are four circles tangent to lines $B'C$, $C'B$, and also tangent to the circumcircles of triangles $A'BC$ and $A'B'C'$.

	Let $P$ be another point of intersection of line $C'B$ and circle $A'B'C'$.  Similarly $Q$ is the other point of intersection of line $CB'$ and circle $BCA'$.  Then using the inscribed angle theorem, we have
	\begin{gather*}
		\angb B'PB=\angb B'A'C'=\angb B'QB \text{ and }
		\angb PB'A'=\angb PC'A'=\angb BB'Q.
	\end{gather*}

	We see that quadrilateral $B'PBQ$ is a deltoid (kite), and two of its sides are chords of circles $A'BC$ and $A'B'C'$. Note also that
	\begin{equation*}
		\angb PC'B'=\angb BCB'.
	\end{equation*}
	That is, the central angles of arcs $PB'$ and $QB$ of circles $A'C'B'$ and $A'BC$ are equal. And hence, by the segment theorem for lines $C'B$ and $CB'$ there are four circles (two for the inscribed circle of the deltoid and two for its escribed circle) that are tangent to the lines $B'C$ and $C'B$ and also to the circles $A'BC$ and $A'B'C'$.
\end{proof}

\section{Monge's theorem and a property of the Feuerbach point}
\label{sct:(monge theorem)}

We show that in hyperbolic geometry Monge's theorem holds, in some sense dual to the theorem on the radical axes of three circles. It is very helpful in many situations. In particular, it reveals a remarkable property of the Feuerbach point.

\begin{theorem}[Three dunce caps theorem]
	Let $\omega_1$, $\omega_2$ and $\omega_3$ be three circles in the hyperbolic plane. Let $P_1$ be the intersection point of the external tangents to circles $\omega_2$ and $\omega_3$. We define points $P_2$ and $P_3$ similarly. Then the points $P_1$, $P_2$ and  $P_3$ are collinear. The theorem is still valid when exactly two of the points are the intersections of the internal tangents.
\end{theorem}

This theorem can be proved directly using the hyperbolic analogue of Menelaus's theorem. However, it is more elegant to use classical Euclidean arguments after conversion to three-dimensional space.

\begin{proof}
	Let's construct on each circle a sphere of the same radius, and consider a plane tangent to all three spheres. It is easy to see that this plane passes through the points $P_1$, $P_2$ and $P_3$, since it contains the generators of the cones. However, the intersection of this plane and the plane of the circles is a straight line on which must lie $P_1$, $P_2$ and $P_3$.
\end{proof}

Since the points $P_1$, $P_2$ and $P_3$ are the centers of homothety for the corresponding circles, we have the following theorem.

\begin{theorem}[Monge's theorem]
	Let $\omega_1$, $\omega_2$ and $\omega_3$ be three circles. Let $P_1$ be the center of a positive homothety from circles $\omega_2$ and $\omega_3$. We define points $P_2$ and $P_3$ similarly. Then the points $\omega_1$, $\omega_2$ and $\omega_3$ are collinear. The theorem is still valid if exactly two of the three homotheties are negative.
\end{theorem}
In the proof, we essentially relied on the fact that the center of this homothety exists (and if nonexistent then we employ analytic continuation). However, for intersecting circles or in the case of one circle inside another, a solution is not, generally speaking, so evident. Let us show how to avoid these problems.

The {\em homothetic center} of the circles $\omega_1$ and $\omega_2$ is designated as that point, such that lines passing through it cross $\omega_1$ and $\omega_2$ in equal angles. There are two such points and, depending on whether the senses of the angles relative to the line are the same or opposite, one center is positive and the other negative. 

It is easy to show that if two lines through a point intersect $\omega_1$ and $\omega_2$ in equal angles, then all lines through the point have this property (it suffices to construct a circle symmetric about the diameter of the Poincar\'e disk). One of these lines is the line passing through the centers of the circles. It crosses $\omega_1$ and $\omega_2$ at right angles. A second such line is the diameter of the Poincar\'e disk passing through the Euclidean homothety center (positive or negative) of the two circles (without loss of generality, we assume that the centers of both circles do not lie on the same diameter). This line crosses the circles at equal angles, and therefore this line in the Poincar\'e disk has the same property in the hyperbolic plane. We have obtained two separate lines intersecting $\omega_1$ and $\omega_2$ in equal angles, so the intersection of these lines is the center of homothety.

This immediately implies Monge's theorem: a line intersecting pairs of circles $\omega_1$ and $\omega_2$ and also $\omega_2$ and $\omega_3$ in equal angles also intersects circles $\omega_1$ and $\omega_3$ at equal angles, i.e., passes through the center of homothety.

Monge's theorem is extremely useful for proofs involving the concurrence of lines or the collinearity of points. Here is one such usage.

\begin{theorem}
	Given triangle $ABC$ and circle $\omega$, let circle $\omega_a$ be tangent to the rays $AB$ and $AC$ and also internally tangent to $\omega$ at point $P_a$. We similarly define circles $\omega_b$ and $\omega_c$ and points $P_b$ and $P_c$. Then the lines $AP_a$, $BP_b$ and $CP_c$ are concurrent.
	\label{thm:no13}
\end{theorem}
\begin{proof}
	Let the point $P$ be the positive homothety center of the circle $\omega$ and the inscribed circle of triangle $ABC$. Then from Monge's theorem, $A$, $P$ and $P_a$ are collinear. We apply a similar argument to $BP_b$ and $CP_c$. Thus the three lines concur at point $P$.
\end{proof}

From this demonstration it is clear that this theorem admits a number of variations. For example, we can require that all of the other circles are externally tangent to $\omega$. Then the point of concurrence will be the negative center of homothety of the inscribed circle and $\omega$.

Since the homothety center of the inscribed circle and the Euler circle lies on the line $EI$, Theorem~\ref{thm:no13} implies the following:

\begin{theorem}
	Given triangle $ABC$ in the hyperbolic plane, let $F$, $F_a$, $F_b$ and $F_c$ be the points of contact of the Euler circle with the inscribed circle and the correspondeds escribed circles. Suppose~$I$ is a incenter of the triangle. Then the lines $FI$, $AF_a$, $BF_b$ and $CF_c$ are concurrent.
\end{theorem}

\section{What happens on the sphere?}

The reader who is familiar with hyperbolic geometry knows that in an ``algebraic sense'' hyperbolic geometry and spherical geometry have the same structure. Therefore everything described in this article can also be applied to spherical geometry. This might be simpler in some sense for two reasons. First, spherical geometry is ``real''. The sphere is easy to imagine or touch (take, for example, a globe). It does not require complicated definitions, and theorems seem more visible and motivated. But this surface has one significant drawback: it is very inconvenient to draw on a plane. It can be done, for example, by a stereographic projection from which we get the ``sphere on a plane'' model, which has many similarities to the Poincar\'e disk. But the resulting model is still inconvenient to work with, at least because lines is hard to distinguish ``to the eye'' from usual circles (lines in this model are circles if the power of the point $(0,0)$ with respect to the circle is~$-1$).

To illustrate the similarity of the geometries, we prove for the sphere Lexell's theorem on the locus of the set of points $X$ such that the area of triangle $ABX$ is constant. First of all, note that equidistant curves on a sphere are ordinary circles (this is where indeed you see that equidistant curves and circles are ``one and the same''). If you look at a globe and consider that the equator is a straight line, the parallels of latitude are equidistant curves, but also circles centered at the north or south poles.

The second question that occurs is what is the ``analogue'' of the points $A^*$ and $B^*$, which are the inverses of $A$ and $B$ in the Poincar\'e disk model? Note that the point $A$ has the property that any line through $A$ (which is a circle in the Poincar\'e disk) also passes through the point $A^*$. On the sphere, this property attaches to the antipode of $A$, the point symmetric to $A$ with respect to the center of the sphere.

We now state the analogue of Theorem~\ref{thm:parallel} on the sphere.

\begin{theorem}[Lexell's theorem]
	Let $A^*$ and $B^*$ be points opposite $A$ and $B$ with respect to the center of the sphere. Then for any point $X$ on any arc $A^*B^*$, the area of triangle $ABX$ is constant, i.e., will not depend on the choice of $X$.
\end{theorem}
\begin{proof}
	There is a very simple and elegant proof that the set of points $X$ consists of an equidistant curve. Here we present a sketch; a detailed proof can be found in \cite{atanasyan2001geometry_en} or \cite{Bibkov2007ontrisection_en}.

	\begin{center}
		\includegraphics{figfeuerbach-5}
		
		\r\label{last}
		\end{center}

	Given triangle $ABX$ (Fig.~\ref{last}), denote by $A_1$ and $B_1$ the midpoints of $AX$ and $BX$ respectively, and by $l$ the line through these points. Let $P$, $Q$ and $R$ be the orthogonal projections of points $A$, $B$ and $X$ on $l$. Note that the right triangles $APA_1$ and $XRA_1$ are congruent because the hypotenuse and one acute angle of each are equal. The same is true of triangles $BQB_1$ and $XRB_1$. From this we can draw two conclusions. First,
	\begin{equation*}
		S_{\trib A_1XB_1} = S_{\trib APA_1}+S_{\trib BQB_1}, 
		\text{ therefore }
		S_{\trib AXB} = S_{{}_\square APQB}\text{.}
	\end{equation*}
	And second, $AP=BQ=XR$, i.e., points $A$ and $B$ lie on a curve equidistant from $l$ and $X$ lies on a symmetric equidistant curve. We denote these equidistant curves by $\omega'$ and $\omega$ respectively.

	We can deduce that for every point $X'$ on the arc $A^*B^*$ the area of triangle $AX'B$ is equal to the area of quadrilateral $APQB$ (note that equality of required right triangles $APA'_1$ and $X'R'A'_1$ holds because its have a leg and an acute angle that are equal).

	We know that on the sphere curves equidistant to any straight line (great circle) are ordinary circles, ``parallel'' to it. It is easy to see that $\omega$ and $\omega'$ are symmetrical not only with respect to $l$, but also with respect to the center of the sphere, so that $\omega$ will pass through the points $A^*$ and $B^*$ (i.e., will be the very circle indicated in the theorem).
\end{proof}

\begin{center}
\includegraphics{figfeuerbach-6}
 
\r
\end{center}

For hyperbolic geometry, this last statement of the proof requires a small change, as it is unclear what the analogue of ``symmetric with respect to the center of the sphere'' is. The interested reader can examine how the reflection in center of the sphere appears in the stereographic projection and discover that it is an inversion in a circle of radius $i$. Therefore in the Poincar\'e disk model, it is logical to take this transformation as an inversion with respect to the absolute (a circle of radius~$1$). Then we see that $\omega$ and $\omega'$ are {\em inverses with respect to the absolute} in the Poincar\'e disk model. This is particularly easy to comprehend if we employ the upper half plane model and let $l$ be a ray perpendicular to the line of the absolute.

\subsection*{\begin{center}Acknowledgements \end{center}}
I thank Fedor Petrov and Mikhail Vyalyi for their interest in this work and useful discussions. I am deeply grateful to Robert~A.~Russell for the translation of this text from Russian, his remarks and revisions making this text clear and correct.


\end{document}